\documentclass{article}

\usepackage{amsmath}
\usepackage{amssymb}
\usepackage{amsthm}
\usepackage{graphicx}
\usepackage{rotating}
\usepackage{bm}

\newtheorem{theorem}{Theorem}
\newtheorem{proposition}[theorem]{Proposition}
\newtheorem{lemma}[theorem]{Lemma}
\newtheorem{definition}[theorem]{Definition}
\newtheorem{example}[theorem]{Example}

\newtheorem{conjecture}[theorem]{Conjecture}

\begin{document}

\markboth{Stefan Gerhold}{Point lattices and oscillating recurrence sequences}

\title{Point lattices and oscillating recurrence sequences}

\author{Stefan Gerhold
 \thanks{Supported by the SFB-grant F1305 of the Austrian FWF}
 \\ RISC, J. Kepler University Linz
 \\ Altenbergerstr. 69, 4040 Linz, Austria
 \\ stefan.gerhold@risc.uni-linz.ac.at
 \\Tel. +43 732 2468 9969}

\maketitle

\begin{abstract}
We consider the following question: Which real sequences $(a(n))$ that satisfy
a linear recurrence with constant coefficients are positive for sufficiently large $n$? We show that
the answer is negative for both $(a(n))$ and $(-a(n))$,
if the dominating characteristic roots in the representation
of $a(n)$ as a generalized power sum comprise either
two pairs of conjugate complex roots or
one negative real root and one pair of conjugate complex roots.
The proof uses results from Diophantine approximation and the geometry of numbers.
Concerning the case of a real positive dominating root we show what the
answer to the question of positivity is for almost all values of the other dominating roots,
provided that all dominating roots are simple.
\end{abstract}



\section{Introduction and main result}

Several classical special function inequalities, such as Fej\'{e}r's Inequality \cite{An00} or the Askey-Gasper Inequality \cite{An00},
assert the positivity of an object that can be defined by a linear recurrence with polynomial
coefficients. Even for the special case of linear recurrences
\begin{equation}\label{eq:cfinite rec}
  a(n+d) = s_1 a(n+d-1) + \dots + s_{d-1}a(n+1) + s_d a(n), \quad n\in\mathbb{N},
\end{equation}
with constant coefficients $s_1,\dots,s_d\in\mathbb{R}$ it is not
always a simple matter to decide from the recurrence coefficients and the real initial values
$a(0),\dots,a(d-1)$ whether the solution is positive or not.
We call sequences $(a(n))$ that satisfy a recurrence of the form \eqref{eq:cfinite rec} {\em recurrence sequences}.
Zeilberger \cite{Ze90} gives them the more suggestive name {\em $C$-finite sequences\/}.
Linear combinations (with constant coefficients) of recurrence sequences are recurrence sequences again,
so positivity results are useful for comparing the magnitude of two sequences, too.

It is well known \cite{Ev03} that 
the sequence $(a(n))$ can be written in terms of the roots $\alpha_1,\dots,\alpha_s$ of the
{\em characteristic polynomial\/}
\[
  z^d - s_1 z^{d-1}- \dots - s_{d-1} z - s_d
\]
of the recurrence as a generalized power sum
\begin{equation}\label{eq:a(n) repr}
  a(n) = C_1(n) \alpha_1^n + \dots C_d(n) \alpha_s^n,
\end{equation}
where the $C_k(n)$ are polynomials in $n$ with complex coefficients. Given a recurrence of the form \eqref{eq:cfinite rec}
and initial values $a(0),\dots,a(d-1)$, the $\alpha_k$ and the $C_k$ can be readily computed.
We refer to the $\alpha_k$ that occur in \eqref{eq:a(n) repr} with nonzero coefficient as {\em characteristic roots of\/}
$(a(n))$. The characteristic roots of maximal modulus will be called {\em dominating characteristic roots of\/}
$(a(n))$.

\begin{example}\label{ex:cfinite rec}
  Consider the recurrence
  \begin{align*}
    a(n+5) & =  3 a(n+4) -2(\sqrt{5}+1)a(n+3)+6(\sqrt{5}+1)a(n+2) \\
           & \phantom{=} -16a(n+1) + 48a(n).
  \end{align*}
  Its characteristic polynomial is
  \[
    (z-3)(z-2 \mathrm{e}^{7\mathrm{i}\pi/5})(z-2\mathrm{e}^{-7\mathrm{i}\pi/5})
           (z-2\mathrm{e}^{2\mathrm{i}\pi/5})(z-2\mathrm{e}^{-2\mathrm{i}\pi/5}),
  \]
  and the solution is given by
  \begin{align}
    a(n) & = c_0 3^n + 2^n \left(c_1 \mathrm{e}^{7n\mathrm{i}\pi/5} + \overline{c}_1 \mathrm{e}^{3n\mathrm{i}\pi/5}
         + c_2 \mathrm{e}^{2n\mathrm{i}\pi/5} + \overline{c}_2 \mathrm{e}^{8n\mathrm{i}\pi/5} \right) \notag \\
      & = c_0 3^n + O(2^n), \label{eq:O(2^n)}
  \end{align}
  where the coefficients $c_0\in\mathbb{R}$ and $c_1,c_2\in\mathbb{C}$ depend on the real
  initial values $a(0),\dots,a(3)$.
  We may ask ourselves whether $3^{-n}a(n)$ approaches $c_0$ from one side only. If $c_1$ and $c_2$
  do not both vanish, it is natural to expect (and will be established in this paper) that this does not hold, because
  the $O(2^n)$ term seems to oscillate.
\end{example}

We pose the following conjecture.

\begin{conjecture}\label{cj:equ abs}
  Let $(a(n))$ be a recurrence sequence with no real positive dominating characteristic root. Then
  there are infinitely many $n$ with $a(n)>0$ and infinitely many $n$ with $a(n)<0$.
\end{conjecture}

The sequence $(a(n))$ might not oscillate if there is a real positive dominating characteristic root.
See Section~\ref{se:pos real} for more on this.

So far Conjecture~\ref{cj:equ abs} has only been verified for one dominating characteristic root (trivial) and
for one pair of conjugate complex roots \cite{Bu81}.
We cannot follow an argument from Nagasaka and Shiue \cite{Na90}, viz.\  that this special case
should immediately imply the truth of the conjecture in general. 
The main goal of this paper
is to establish the following theorem by an extension of Burke and Webb's proof.

\begin{theorem}[Main Theorem]\label{thm:main}
  Let $(a(n))$ be a recurrence sequence with at most four dominating characteristic roots, none of which is real positive. Then
  there are infinitely many $n$ with $a(n)>0$ and infinitely many $n$ with $a(n)<0$.
\end{theorem}

The rest of the paper is organized as follows. In Section~\ref{se:pre} we reduce Theorem~\ref{thm:main}
from multiple roots to simple roots and subsequently
to a geometric statement about the distribution modulo one of integer multiples of a real vector
$\bm{\xi}=(\xi_1,\xi_2)=(\arg{\alpha_1},\arg{\alpha_2})/2\pi$,
except for some special cases of Theorem~\ref{thm:main} that are settled in Section~\ref{se:compl}.
In Section~\ref{se:irr} we deduce the
desired result from Kronecker's Approximation Theorem, provided that one of $\xi_1,\xi_2$ is irrational. The proof in the case
where both are rational is the subject of Sections~\ref{se:rat} and \ref{se:edrat}.
Section~\ref{se:pos real} presents a metric result that deals with the case of a positive real root.
In the conclusion we comment, among other things,
on extending our approach to Conjecture~\ref{cj:equ abs} to an arbitrary number of dominating characteristic roots.

\section{Notation and preliminaries}\label{se:pre}

We write $\mathbb{N}$, $\mathbb{Z}$, $\mathbb{Q}$, $\mathbb{R}_0^+$, $\mathbb{R}$ and $\mathbb{C}$
for the sets of natural numbers (including zero), integers, rational numbers, non-negative real numbers,
real numbers and complex numbers, respectively. The conjugate of a complex number $z$ is denoted by $\overline{z}$.
Whenever $\bm{v}$ is a vector, we use the same letter with a subscript for its components,
as in $\bm{v}=(v_1,\dots,v_m)$.
For vectors $\bm{\xi}=(\xi_1,\dots,\xi_m)$, $\bm{\rho}=(\rho_1,\dots,\rho_m)$ of real numbers
and a real number $\rho$, we write
\begin{align*}
  \bm{\xi} \bmod \bm{\rho} &= (\xi_1 \bmod \rho_1,\dots, \xi_1 \bmod \rho_1) \quad \textnormal{and} \\
  \bm{\xi} \bmod \rho &= (\xi_1 \bmod \rho,\dots, \xi_1 \bmod \rho).
\end{align*}
We define the open rectangle parallel to the axes with side lengths $2\lambda_1,2\lambda_2\in\mathbb{R}$
centred at $\bm{c}=(c_1,c_2)\in\mathbb{R}^2$ as
\[
  \mathcal{R}_{\lambda_1,\lambda_2}(\bm{c}) := \left\{ \bm{x} \in \mathbb{R}^2 :
  |x_1-c_1|<\lambda_1, |x_2-c_2|<\lambda_2 \right\}.
\]
For an open square parallel to the axes we write
\[
  \mathcal{S}_\lambda(\bm{c}) := \mathcal{R}_{\lambda,\lambda}(\bm{c}),\quad \lambda \in \mathbb{R}, \bm{c} \in \mathbb{R}^2.
\]
A {\em lattice} \cite{Ca59} is a discrete subgroup $\Lambda\subset\mathbb{R}^m$. Its determinant
is denoted by $\bm{d}(\Lambda)$.
The lattice $L_g(a_1,a_2)$ is defined in Section~\ref{se:rat}.
\vskip 5mm

Let $(a(n))$ be as in Theorem~\ref{thm:main}.
We order the characteristic roots $\alpha_1,\dots,\alpha_{s}$ of $(a(n))$ such that $\alpha_1,\dots,\alpha_{l}$ contain
all real dominating characteristic roots, precisely one element of every pair of conjugate non-real dominating characteristic roots
and no other roots.
Note that this implies $l=1$ or $l=2$.

Moreover, let $\alpha_1,\dots,\alpha_{l}$ be ordered such that 
\[
  D:= \deg C_1 = \dots = \deg C_m > \deg C_{m+1} \geq \dots \geq \deg C_{l}
\]
for some $1\leq m \leq l\leq 2$. Then we obtain \cite{Kn94}
\[
  n^{-D}a(n) = \sum_{k=1}^m\left(c_k\alpha_k^n +
  \overline{c}_k\overline{\alpha}_k^n\right) + O(n^{-1}|\alpha_1|^n),
\]
where $c_k$ is the leading coefficient of $C_k(n)$.
This formula shows that Theorem~\ref{thm:main} can be deduced from Burke and Webb's result
($m=1$) and the following theorem ($m=2$). Observe that we can safely assume $|\alpha_1|=|\alpha_2|=1$,
since we can divide by the positive factor $|\alpha_1|^n$.

\begin{theorem}\label{thm:main2}
  Let $\alpha_1,\alpha_2\in\mathbb{C}\backslash\mathbb{R}_0^+$, $|\alpha_1|=|\alpha_2|=1$,
  $\alpha_1\neq\alpha_2\neq\overline{\alpha}_1$. Let further $c_1,c_2$ be nonzero complex numbers and
  \begin{equation}\label{eq:b(n)}
    b(n) := c_1 \alpha_1^n + \overline{c}_1\overline{\alpha}_1^n
           + c_2 \alpha_2^n + \overline{c}_2\overline{\alpha}_2^n, \qquad n\geq 0.
  \end{equation}
  Then there is $\delta>0$ such that $b(n)>\delta$ for infinitely many $n$ and
  $b(n)<-\delta$ for infinitely many $n$.
\end{theorem}

Note that if $\delta$ was replaced by zero, it might happen that
e.g.\ all negative values $b(n)$ are so small in absolute value that the remainder term of $a(n)$, which comes from
the characteristic roots of smaller modulus, takes over and makes the corresponding values
$a(n)$ positive.
This uniformity condition was missed by Burke and Webb \cite{Bu81}.
They only argue that $c_1 \alpha_1^n + \overline{c}_1\overline{\alpha}_1^n$ has
infinitely many positive and infinitely many negative values, which is not sufficient, but their proof can be easily repaired.

Now let $\alpha_1$, $\alpha_2$, $c_1$, $c_2$ be as in Theorem~\ref{thm:main2}.
Replacing $(\alpha_k,c_k)$ by $(\overline{\alpha}_k,\overline{c}_k)$ and vice versa if necessary, we may assume $\mathrm{Im}(c_k)\geq 0$.
Putting $\theta_k:=\arg{\alpha_k}$, we obtain by standard formulas
\begin{align*}
  b(n) & = 2\sum_{k=1}^2 \mathrm{Re}\left( c_k \exp\left(\mathrm{i}n\theta_k \right) \right) \\
  & =  2\sum_{k=1}^2\left(\mathrm{Re}(c_k) \cos{n\theta_k} - \mathrm{Im}(c_k) \sin{n \theta_k}\right) \\
  & = \sum_{k=1}^2 w_k \sin(n\theta_k + \varphi_k),
\end{align*}
where the coefficients are nonzero real numbers
\[
  w_k :=
  \begin{cases}
    -2|c_k|, & c_k \in \mathbb{C}\backslash\mathbb{R}; \\
    2c_k, & c_k \in \mathbb{R},
  \end{cases}
\]
and the $\varphi_k$ are given by
\[
  \varphi_k :=
 \begin{cases}
  -\arctan\tfrac{\mathrm{Re}(c_k)}{\mathrm{Im}(c_k)}, & c_k \in \mathbb{C}\backslash\mathbb{R}; \\
  \frac{1}{2}\pi, & c_k \in \mathbb{R}.
  \end{cases}
\]
We turn our attention to the signs of $\sin(n\theta_k + \varphi_k)$. If we can prove that
for every pair $(S_1,S_2)$ of $+1$'s and $-1$'s there are
infinitely many $n$ such that the sign of $\sin(n\theta_k + \varphi_k)$
equals $S_k$ for $k=1,2$, we will have shown that $(b(n))$ oscillates,
whatever the values of the $c_k$ (and thus the $w_k$) are. In other words, we are looking for $n$ such that
\[
  (n\theta_k + \varphi_k) \bmod 2\pi \in\ ]0,\pi[
\]
or
\[
  (n\theta_k + \varphi_k) \bmod 2\pi \in\ ]\pi,2\pi[,
\]
respectively. To get the $\delta$ in Theorem~\ref{thm:main2}, we have to shrink the intervals to
\[
  ]\epsilon,\pi-\epsilon[ \quad \text{and} \quad ]\pi+\epsilon,2\pi-\epsilon[
\]
for some small $\epsilon>0$, of course independent from $n$.
Now we rescale to the unit interval.

\begin{theorem}\label{thm:xi0,xi1}
  Let $\xi_1,\xi_2 \in \left]0,1\right[\backslash\{\tfrac{1}{2}\}$ such that $\xi_1\not\equiv\pm\xi_2\pmod{1}$ and,
  if both $\xi_1$ and $\xi_2$ are rational, then the pair of their denominators
  (written with the larger denominator first) is none of
  $(5,5)$, $(6,3)$, $(8,4)$.
  Then for all $\bm{c}\in\mathbb{R}^2$ there is $\epsilon>0$ such that
  there are infinitely many $n$ with
  \[
    n(\xi_1,\xi_2)\bmod 1\ \in\ \mathcal{S}_{1/4-\epsilon}(\bm{c}) \bmod 1.
  \]
\end{theorem}

Since the sine function is continuous, applying this theorem with
$(\xi_1,\xi_2)=(\theta_1/2\pi,\theta_2/2\pi)$ and $c_k=\tfrac{1}{4}-\varphi_k/2\pi$ to make
$\sin(n\theta_k+\varphi_k)$ positive and $c_k=\tfrac{3}{4}-\varphi_k/2\pi$ for a negative sign proves Theorem~\ref{thm:main2},
unless one of the $\alpha_k$ is a negative real number (which implies $\xi_k=\tfrac{1}{2}$)
or $\theta_1/2\pi,\theta_2/2\pi$ are rational numbers with denominators in $\{(5,5), (6,3), (8,4)\}$.
Section~\ref{se:compl} deals with these special cases of Theorem~\ref{thm:main2}.

In the proof of Theorem~\ref{thm:xi0,xi1} we distinguish the following three cases:
\begin{enumerate}
  \item[(1)] $\xi_1, \xi_2, 1$ are linearly independent over $\mathbb{Q}$.  \label{case1}
  \item[(2)] $\xi_1, \xi_2$ are not both rational, but satisfy a linear relation $u_1\xi_1 + u_2\xi_2 = v$ with
  $u_1,u_2,v\in\mathbb{Z}$.  \label{case2}
  \item[(3)] $\xi_1$ and $\xi_2$ are both rational.  \label{case3}
\end{enumerate}

Section~\ref{se:irr} settles the first two cases.
The proof of Theorem~\ref{thm:xi0,xi1} in Case~3 is the content of Sections~\ref{se:rat} and \ref{se:edrat}.

We remark that in order to prove Conjecture~\ref{cj:equ abs} for one
pair of conjugate complex dominating roots,
it suffices to show that for every real number $\xi\neq \tfrac{1}{2}$ with
$0<\xi<1$ and every real number $c$ there is $\epsilon>0$ such that for infinitely many $n$
\[
  n\xi\bmod 1\ \in\ \left]c-\tfrac{1}{4}+\epsilon,c+\tfrac{1}{4}-\epsilon\right[ \bmod 1.
\]
This is essentially what was done (without $\epsilon$, cf. the introduction) by Burke and Webb \cite{Bu81}.

\section{The irrational cases}\label{se:irr}

The closure of the set of integer multiples of a vector $\bm{\xi}=(\xi_1,\xi_2)$ modulo one
is described by a classical result from Diophantine approximation.

\begin{theorem}[Kronecker's Theorem]\label{thm:kronecker}
  Let $\xi_1$, $\xi_2$ be real numbers.
  \begin{itemize}
    \item[(i)] If $\xi_1,\xi_2,1$ are linearly independent over the rationals, then the points
    $n\bm{\xi}\bmod 1$, $n\in\mathbb{N}$, lie dense in the unit square.
    \item[(ii)] If $\xi_1,\xi_2$ are not both rational, but satisfy a relation
    $u_1\xi_1 + u_2\xi_2 = v$ with $u_1,u_2,v\in\mathbb{Z}$ and $\gcd(u_1,u_2,v)=1$, then the points
    $n\bm{\xi}\bmod 1$, $n\in\mathbb{N}$, lie dense on the portions of the lines
    \[
      \ell_t:=\left\{ \bm{x} \in \mathbb{R}^2 : u_1x_1+u_2x_2 = t \right\}, \quad t\in\mathbb{Z},
    \]
    which lie within the unit square.
  \end{itemize}
\end{theorem}
\begin{proof}
  See e.g.\ Niven \cite[Theorems~3.4 and 3.6]{Ni63}.
\end{proof}

\begin{figure}[h]
  \centering
  \includegraphics[scale=.37]{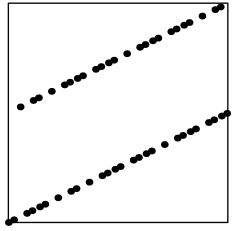}
  \caption{The unit square with $(n\xi_1,n\xi_2)\bmod 1$ for $\xi_1=2\sqrt{2}$, $\xi_1-2\xi_2=2$ and $n=0,\dots,50.$}
  \label{fi:one_rel}
\end{figure}

Part (i) of Theorem~\ref{thm:kronecker} settles Case~1 of Theorem~\ref{thm:xi0,xi1}.
We proceed to Case~2. Let
$\bm{c}\in\mathbb{R}^2$ be arbitrary but fixed and
$\ell_t$ be as in part (ii) of Theorem~\ref{thm:kronecker}.
Since
\[
  \bigcup_{t\in\mathbb{Z}}\ell_t + \mathbb{Z}^2 = \bigcup_{t\in\mathbb{Z}}\ell_t,
\]
it suffices to find infinitely many $n\bm{\xi}\bmod 1$ in the set
\[
  \mathcal{S}_{1/4-\epsilon}(\bm{c}) \cap \bigcup_{t\in\mathbb{Z}} \ell_t,
\]
where $\epsilon>0$ is yet to be chosen.
First suppose that $\xi_1$ and $\xi_2$ are irrational.
Then the parallel lines $\ell_t$ are neither horizontal nor vertical,
since $u_1u_2\neq 0$. Two adjacent lines $\ell_t$, $\ell_{t+1}$
have horizontal distance $1/|u_1|$
and vertical distance $1/|u_2|$. Since $\xi_1\not\equiv\pm\xi_2\pmod{1}$, one of
these quantities must be smaller than or equal to $\tfrac{1}{2}$. 
Thus
\[
  \mathcal{S}_{1/4}(\bm{c}) \cap \bigcup_{t\in\mathbb{Z}} \ell_t \neq \emptyset.
\]
In fact this set is not only non-empty but contains a line segment.
Clearly, we can find $\epsilon>0$ such that the set
$\mathcal{S}_{1/4-\epsilon}(\bm{c}) \cap \bigcup_{t\in\mathbb{Z}} \ell_t$
still contains a line segment of length greater than zero.
Filling this line segment densely with points $n\bm{\xi}\bmod 1$ requires infinitely many $n$.

Now let $\xi_1$ be rational and $\xi_2$ be irrational, and let $b_1\in\mathbb{N}$ be the
denominator of $\xi_1$. This implies $u_2=0$.
Then the lines $\ell_t$ are vertical, and the horizontal distance between $\ell_t$ and $\ell_{t+1}$ is $1/b_1\leq \tfrac{1}{3}$, since
$b_1>2$ by the assumptions of Theorem~\ref{thm:xi0,xi1}. Case~2 of Theorem~\ref{thm:xi0,xi1} is proved.

\section{The rational case}\label{se:rat}

The main goal of this section and the next one is to prove the following theorem.

\begin{theorem}\label{thm:genrat}
  Let $a_1,a_2,b_1,b_2\in\mathbb{N}$, $2\leq b_2\leq b_1$, $1\leq a_k < b_k$,
  $\gcd(a_k,b_k)=1$ for $k=1,2$ and $\tfrac{a_1}{b_1} \not\equiv \pm\tfrac{a_2}{b_2} \pmod{1}$.
  Then there is $\bm{c}\in[0,1]^2$ such that for all $n\in\mathbb{N}$
  \[
    n(\tfrac{a_1}{b_1},\tfrac{a_2}{b_2}) \bmod 1 \notin  \mathcal{S}_{1/4}(\bm{c}) \bmod 1
  \]
  provided that
  \begin{equation}\label{eq:genrat excl}
    (b_1,b_2) \in \left\{ (5,5),(6,3),(8,4) \right\} \cup \left\{ (b_1,2) : 2\leq b_1 \in \mathbb{N} \right\},
  \end{equation}
  and there is no such $\bm{c}$ if \eqref{eq:genrat excl} does not hold.
\end{theorem}

To see that Case~3 of Theorem~\ref{thm:xi0,xi1} follows from Theorem~\ref{thm:genrat}, note that
the purely periodic sequence
\[
  n(\tfrac{a_1}{b_1},\tfrac{a_2}{b_2}) \bmod 1 =
  (\tfrac{na_1\bmod b_1}{b_1},\tfrac{na_2\bmod b_2}{b_2}), \quad n \geq 0,
\]
assumes each of its finitely many values infinitely often.
The $\epsilon$ has disappeared because the set of all $n(\tfrac{a_1}{b_1},\tfrac{a_2}{b_2})\bmod 1$
is finite and $\mathcal{S}_{1/4}(\bm{c})$ is open.

\begin{figure}[h]
  \centering
  \includegraphics[scale=.37]{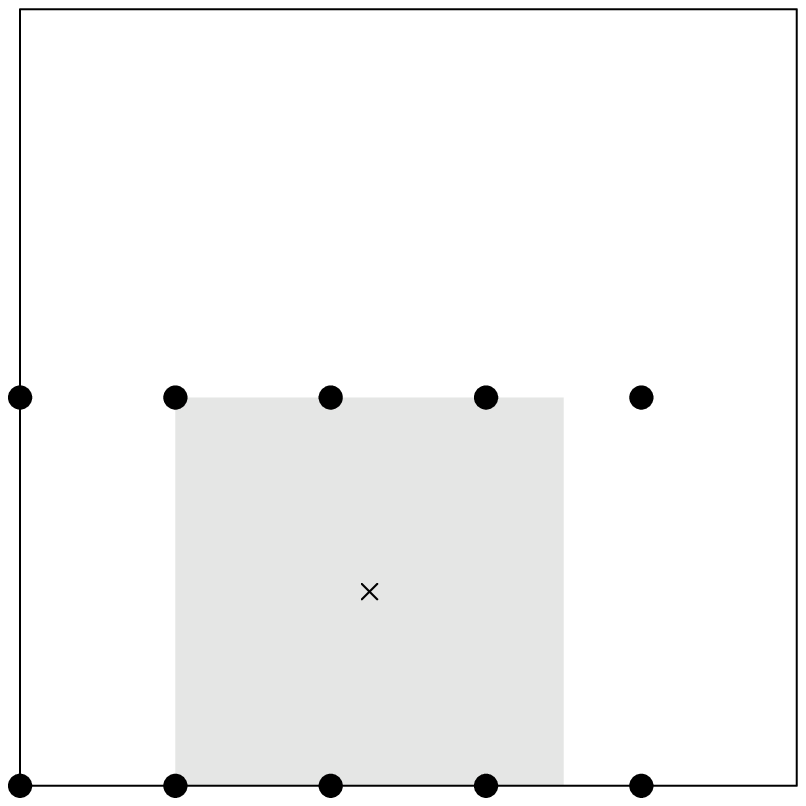} 
  \includegraphics[scale=.37]{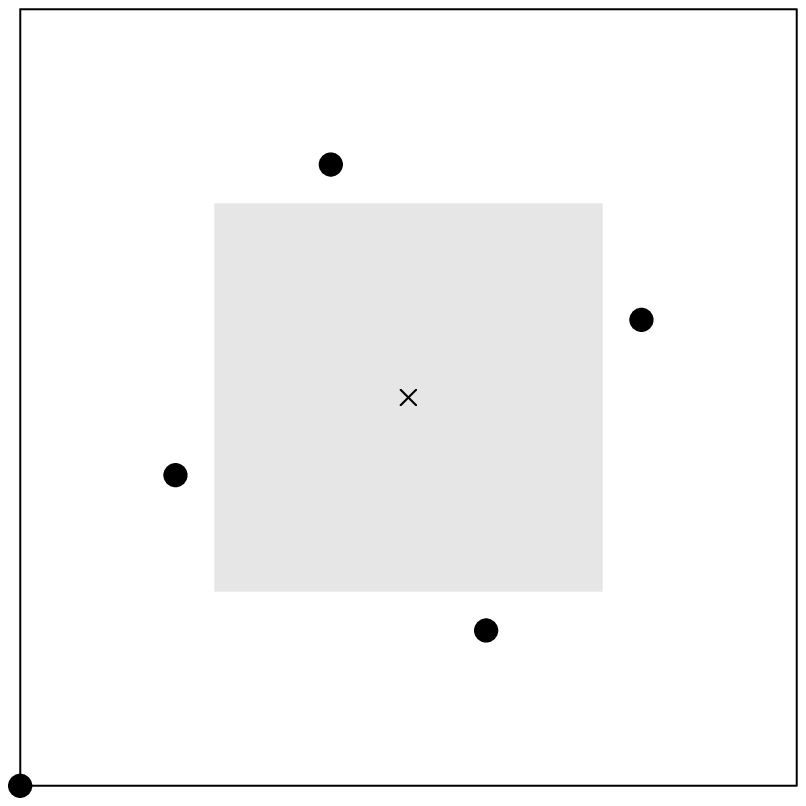}
  \includegraphics[scale=.37]{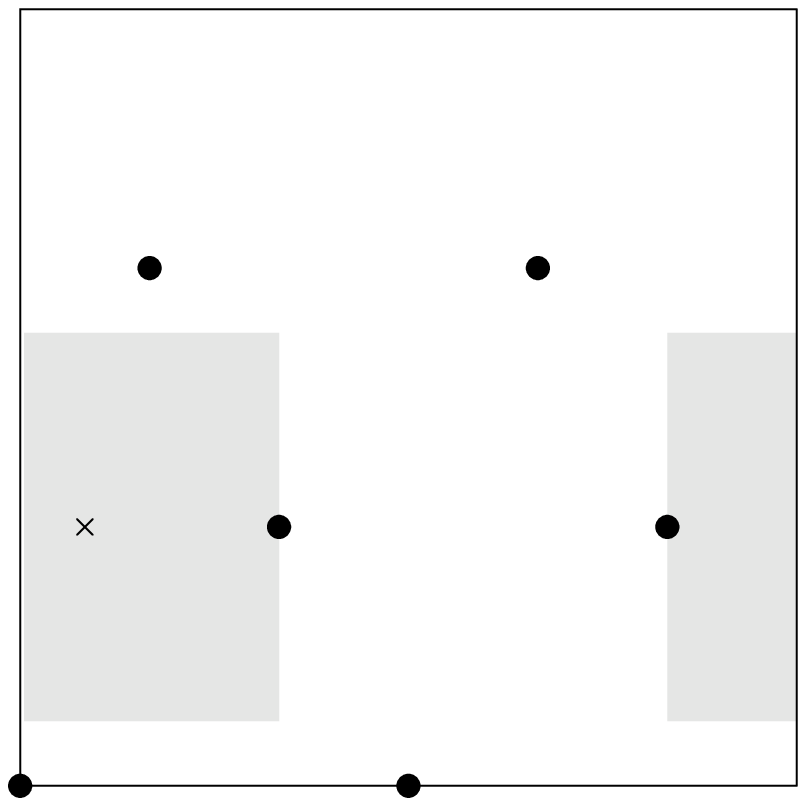}
  \includegraphics[scale=.37]{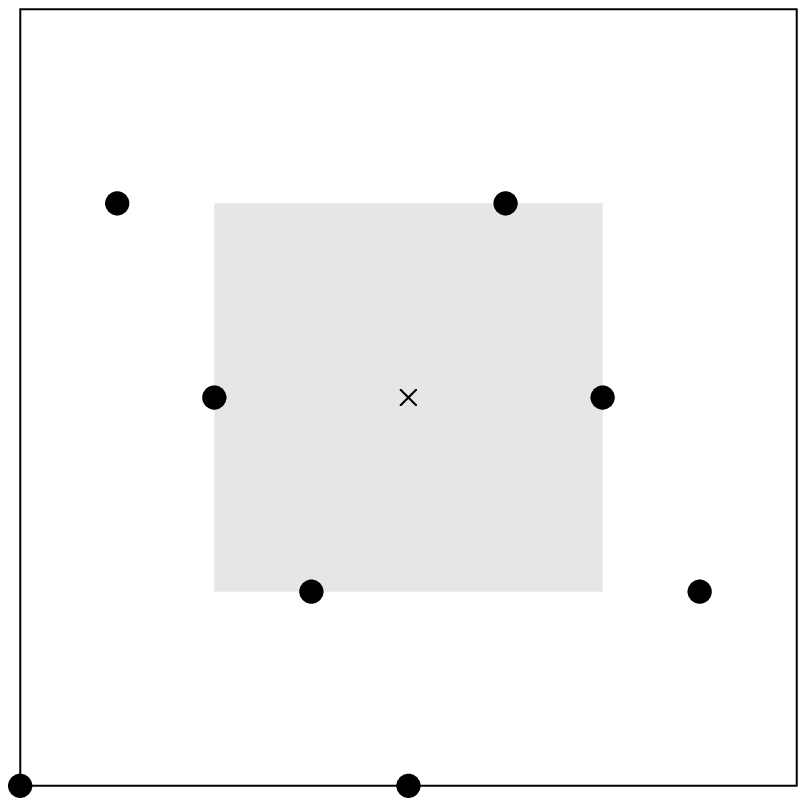}
  \caption{The unit square with the set $\{n(\tfrac{a_1}{b_1},\tfrac{a_2}{b_2})\bmod 1:n\in\mathbb{N}\}$
   for $(\tfrac{a_1}{b_1},\tfrac{a_2}{b_2})=(\tfrac{1}{5},\tfrac{1}{2})$, $(\tfrac{1}{5},\tfrac{2}{5})$,
   $(\tfrac{1}{6},\tfrac{2}{3})$ and $(\tfrac{3}{8},\tfrac{1}{4})$, respectively.
  }\label{fi:spec cases}
\end{figure}

\begin{proof}[Proof of the right to left implication of Theorem~\ref{thm:genrat}]
If $b_2=2$, we necessarily have $a_2=1$, and we may take $c_2=\tfrac{1}{4}$ and $c_1\in\mathbb{R}$ arbitrary.
(See Figure~\ref{fi:spec cases} for an example.)
If $(b_1,b_2)=(5,5)$, it is easy to see that for all $\bm{a}$ in question the set of integer multiples modulo one
is one of the two sets
\[
  \{  n(\tfrac{1}{5},\tfrac{2}{5}) \bmod 1 : n\in\mathbb{N}\} \quad \textnormal{and} \quad
  \{  n(-\tfrac{1}{5},\tfrac{2}{5}) \bmod 1 : n\in\mathbb{N}\},
\]
obtained from $\bm{a}=(1,2)$ and $\bm{a}=(-1,2)$, respectively.
Similarly, for $(b_1,b_2)=(6,3)$ it suffices to consider $\bm{a}=(\pm 1,2)$.
This is also true for $(b_1,b_2)=(8,4)$, if we take $\bm{a}=(\pm 3,1)$ instead
of $(\pm 1,2)$.
The number of $\bm{a}$'s to check can be reduced further
by taking advantage of some obvious symmetries. By the subsequent lemma,
the alternative with negative first entry can be discarded in each of the three cases.
Figure~\ref{fi:spec cases} shows that in the remaining cases we may take $\bm{c}=(\tfrac{1}{2},\tfrac{1}{2})$,
$(\tfrac{1}{12},\tfrac{1}{3})$ and $(\tfrac{1}{2},\tfrac{1}{2})$, respectively.
\end{proof}

\begin{lemma}\label{le:symm}
  Define the maps $s$ and $\tau$ on $\mathbb{R}^2$ by
  \[
    s(x_1,x_2) = ((1-x_1)\bmod 1,x_2) \quad \text{and} \quad \tau(x_1,x_2) = (x_2,x_1).
  \]
  Then for all real numbers $\xi_1$, $\xi_2$
  \[
    s((\xi_1,\xi_2)\bmod 1) = s(\xi_1,\xi_2)\bmod 1 \quad \text{and} \quad \tau((\xi_1,\xi_2)\bmod 1) = \tau(\xi_1,\xi_2)\bmod 1.
  \]
\end{lemma}
\begin{proof}
  Obvious.
\end{proof}

We have shown this implication just for the sake of completeness. The interesting
part of Theorem~\ref{thm:genrat} for our purpose is the converse implication. Its
proof is the content of the remainder of this section and of the following section.

\begin{definition}\label{def:L}
  Let $g$ be a positive integer and $a_1,a_2$ be integers relatively prime to $g$.
  Then we define the {\em lattice of multiples of $\bm{a}=(a_1,a_2)$ modulo $g$\/} as
  \begin{equation*}
    L_g(a_1,a_2) := \left\{ \bm{u}\in\mathbb{Z}^2 : n\bm{a} \equiv \bm{u}\pmod{g}\quad
     \textnormal{for some}\ n \in \mathbb{N}\right\}.
  \end{equation*}
\end{definition}

Alternatively \cite{Ro97}, $L_g(a_1,a_2)$ can be defined as the lattice generated by the vectors
$(0,g),(g,0)$ and $(a_1,a_2)$.

The lattices $L_g(a_1,a_2)$ will provide a convenient representation of the sets of integer
multiples of rational numbers modulo one, which we encountered in Theorem~\ref{thm:genrat}. For this purpose
we require
a version of the well-known Chinese Remainder Theorem for moduli that are not necessarily pairwise relatively prime.

\begin{theorem}[Generalized Chinese Remainder Theorem]\label{thm:gcrt}
  Let $b_1,\dots,b_m$ be positive integers and $u_1,\dots,u_m$ be integers. Then there is an
  integer $0\leq u<\textnormal{lcm}(b_1,\dots,b_m)$ with
  \[
    u \equiv u_i \mod{b_i},\quad 1\leq i\leq m,
  \]
  provided that
  \[
    u_i \equiv u_j \mod{\gcd(b_i,b_j)},\quad 1\leq i,j \leq m.
  \]
\end{theorem}
\begin{proof}
  See Knuth \cite[Exercise~4.3.2.3]{Kn98}.
\end{proof}

\begin{lemma}\label{le:M to L}
  Let $a_1,a_2$ be integers and $b_1,b_2$ be positive integers with $\gcd(a_k,b_k)=1$ for $k=1,2$
  and $g:=\gcd(b_1,b_2)$.
  Then
  \begin{flalign*}
     (i)& \quad \left\{ n(\tfrac{a_1}{b_1}, \tfrac{a_2}{b_2})\bmod 1 : n\in\mathbb{N}\right\}
      = \left\{ (\tfrac{u_1}{b_1},\tfrac{u_2}{b_2}) : \bm{u}\in L_g(a_1,a_2), 0\leq u_k < b_k\right\} & \\
     (ii)& \quad L_g(a_1,a_2) = \left\{ \bm{u} \in \mathbb{Z}^2 : a_1u_2 \equiv a_2u_1\pmod{g} \right\} &
  \end{flalign*}
\end{lemma}
\begin{proof} We have
  \begin{align*}
    &\phantom{=} \left\{ n(\tfrac{a_1}{b_1}, \tfrac{a_2}{b_2})\bmod 1 : n\in\mathbb{N}\right\} \\
    &= \left\{ (\tfrac{n a_1 \bmod b_1}{b_1}, \tfrac{n a_2\bmod b_2}{b_2}) : n\in\mathbb{N}\right\} \\
    &= \left\{ (\tfrac{u_1}{b_1},\tfrac{u_2}{b_2}) : n\bm{a}\equiv \bm{u}\pmod{\bm{b}},
     0\leq u_k<b_k,\ k=1,2,\quad \textnormal{for some}\ n \in \mathbb{N}\right\} \\
    &= \left\{ (\tfrac{u_1}{b_1},\tfrac{u_2}{b_2}) : \bm{u}\in L_g(a_1,a_2), 0\leq u_k < b_k\right\}.
  \end{align*}
  The latter equality and assertion (ii) follow from Theorem~\ref{thm:gcrt}.
\end{proof}

\begin{figure}[h]
  \centering
  \includegraphics[scale=.6]{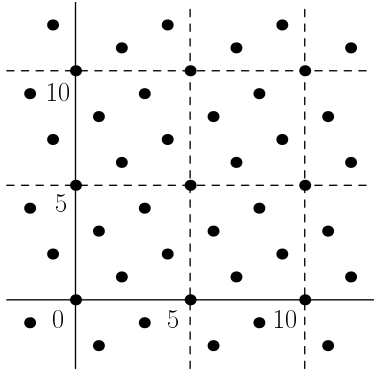}
  \caption{The lattice $L_5(2,1)$.}
  \label{fi:example L}
\end{figure}

\begin{example}\label{ex:example L}
  In Example~\ref{ex:cfinite rec} the $O(2^n)$ term yields $\xi_1=\theta_1/2\pi=\tfrac{7}{10}$
  and $\xi_2=\theta_2/2\pi=\tfrac{2}{5}$. The corresponding lattice $L_5(7,1)=L_5(2,1)$ is displayed
  in Figure~\ref{fi:example L}.
\end{example}

Let $\bm{a}=(a_1,a_2)$ and $\bm{b}=(b_1,b_2)$ be as in the assumptions of Theorem~\ref{thm:genrat},
but such that $\bm{b}$ is not in the set~\eqref{eq:genrat excl}, and put $g:=\gcd(b_1,b_2)$.
In the light of Lemma~\ref{le:M to L}, it is an immediate consequence of the periodicity property
\begin{equation}\label{eq:period}
  L_g(a_1,a_2) = L_g(a_1,a_2) + g\mathbb{Z}^2
\end{equation}
that searching a point $n(\tfrac{a_1}{b_1},\tfrac{a_2}{b_2})\bmod 1$ in a `modded' square $\mathcal{S}_{1/4}(\bm{c})\bmod 1$
amounts to looking for a point of the lattice $L_g(a_1,a_2)$ in the rectangle
$\mathcal{R}_{b_1/4,b_2/4}(b_1 c_1,b_2 c_2)$ with side lengths $b_1/2$, $b_2/2$.
We let $c_k$ absorb $b_k$ and write again $\bm{c}=(c_1,c_2)$ for the arbitrary centre $(b_1 c_1,b_2 c_2)$.

\begin{example}\label{ex:point in rect}
  If we want to show that the $O(2^n)$ term in \eqref{eq:O(2^n)} oscillates, we are lead to the problem
  of finding a point of $L_5(2,1)$ in any rectangle $\mathcal{R}_{5/2,5/4}(\bm{c})$, $\bm{c}\in\mathbb{R}^2$.
\end{example}

If the numbers $b_1/g$ and $b_2/g$ are large, it is easy to find a point
of $L_g(a_1,a_2)$ in the rectangle, whereas $b_1=b_2=g$ is the most difficult case.
This is so because if we fix $a_1$, $a_2$ and $g$ and enlarge $b_1/g$ and $b_2/g$, the lattice
$L_g(a_1,a_2)$ remains invariant, while the rectangle becomes bigger.

At first glance, the problem seems to be easily reducible to the case of equal denominators $b_1=b_2=g$.
In Example~\ref{ex:point in rect}, if we could show that any square $\mathcal{S}_{5/4}(\bm{c})$ contains a point
of $L_5(2,1)$, then it would follow at once that every rectangle $\mathcal{R}_{5/2,5/4}(\bm{c})$ contains
a point of $L_5(2,1)$. But we have already seen (Theorem~\ref{thm:genrat}) that there are squares $\mathcal{S}_{5/4}(\bm{c})$
without points of $L_5(2,1)$. In general, the catch is that
even if $(a_1,a_2,b_1,b_2)$ satisfy the requirements of Theorem~\ref{thm:genrat}
and $(b_1,b_2)$ is not in the set \eqref{eq:genrat excl},
it may happen that $(a_1\bmod g,a_2\bmod g,g,g)$
violate the requirements of Theorem~\ref{thm:genrat} or that $(g,g)$ is in \eqref{eq:genrat excl}.
Therefore we choose a different approach for the case $b_1\neq b_2$.

For relatively prime $b_1$ and $b_2$ the lattice $L_g(a_1,a_2)$ equals $\mathbb{Z}^2$.
All rectangles $\mathcal{R}_{b_1/4,b_2/4}(\bm{c})$ with $\bm{c}\in\mathbb{R}^2$
have side lengths greater than one and therefore contain a point of $\mathbb{Z}^2$.
If $g=2$, then $a_1$ and $a_2$ must be odd,
hence
\[
  L_g(a_1,a_2) = \left\{ \bm{u} \in \mathbb{Z}^2 : u_1\equiv u_1\pmod{2} \right\}.
\]
Since $b_1>4$ in this case, it is easy to see that this lattice contains a point of any rectangle
$\mathcal{R}_{b_1/4,b_2/4}(\bm{c})$.

From now on we assume $g \geq 3$.
The following proposition deals with the case $(b_1,b_2)=(2g,g)$.
Recall that $(b_1,b_2)=(4,2)$, $(6,3)$ and $(8,4)$ need not be considered, because
they are in the set~\eqref{eq:genrat excl}.

\begin{proposition}\label{pr:2g 2g}
  Let $a_1,a_2,b_1,b_2$ be as in Theorem~\ref{thm:genrat}.
  Suppose $g\geq 5$, $b_1=2g$ and $b_1=g$. Then for all $\bm{c}\in\mathbb{R}^2$
  \[
    L_g(a_1,a_2) \cap \mathcal{R}_{b_1/4,b_2/4}(\bm{c}) \neq \emptyset.
  \]
\end{proposition}

\begin{proof}
  Observe that by the periodicity property \eqref{eq:period} of $L_g(a_1,a_2)$
  it suffices to find a point of the lattice in the set
  \begin{equation}\label{eq:3g set}
    \mathcal{R}_{b_1/4,b_2/4}(\bm{c}) + g\mathbb{Z}^2.
  \end{equation}
  Let $\bm{p}$ be the lower left corner of $\mathcal{R}_{b_1/4,b_2/4}(\bm{c})$.
  We assume w.l.o.g.\ $0\leq p_1,p_2 < g$
  and define $I:=\ ]p_2,p_2+\tfrac{1}{2}g[$.
  Then \eqref{eq:3g set} contains the set
  \begin{equation}\label{eq:3g strip}
    \left([0,g[ \backslash \{p_1\}\right) \times I = \left([0,g[ \times I\right) \backslash \left(\{p_1\} \times I\right).
  \end{equation}
  The interval $I$ contains at least two integers, since its length is $\tfrac{1}{2}g>2$. Since $a_2$ is invertible
  modulo $g$, there are at least two points of
  $L_g(a_1,a_2)$ in $[0,g[\times I$ by part (ii) of Lemma~\ref{le:M to L}, and
  at least one of them lies in \eqref{eq:3g strip}. 
\end{proof}

Now we consider values of $b_1$ that are at least $3g$, which completes the case $b_1\neq b_2$
of Theorem~\ref{thm:genrat}.

\begin{proposition}\label{pr:3g}
  Let $a_1,a_2,b_1,b_2$ be as in Theorem~\ref{thm:genrat}.
  Suppose $g\geq 3$ and $b_1\geq 3g$.
  Then for all $\bm{c}\in\mathbb{R}^2$
  \[
    L_g(a_1,a_2) \cap \mathcal{R}_{b_1/4,b_2/4}(\bm{c}) \neq \emptyset.
  \]
\end{proposition}

\begin{proof}
  It suffices to consider $b_1=3g$ and $b_2=g$.
  Proceeding analogously to the proof of Proposition~\ref{pr:2g 2g},
  we arrive at the set $[0,g[ \times I$ instead of \eqref{eq:3g strip}.
  The result follows from part (ii) of Lemma~\ref{le:M to L} and $\tfrac{1}{2}g>1$.
\end{proof}

\section{The rational case with equal denominators}\label{se:edrat}

In order to finish the proof of Theorem~\ref{thm:genrat}, and thus the proof of Theorem~\ref{thm:xi0,xi1},
we will establish the following proposition.

\begin{proposition}\label{pr:g}
  Let $a_1,a_2,b_1,b_2$ be as in Theorem~\ref{thm:genrat}.
  Suppose $b_1=b_2=g\neq 5$.
  Then for all $\bm{c}\in\mathbb{R}^2$
  \[
    L_g(a_1,a_2) \cap \mathcal{S}_{g/4}(\bm{c}) \neq \emptyset.
  \]
\end{proposition}

If $L_g(a_1,a_2)$ contains one or two sufficiently short vectors, its points are dense enough
so that the square $\mathcal{S}_{g/4}(\bm{c})$ is populated by at least one lattice point.
This is the basic idea of our proof of Proposition~\ref{pr:g}.
Although there are algorithms \cite{Le94,Ro97} tailored to $L_g(a_1,a_2)$ for computing a reduced lattice basis,
we do not know of any specialized a priori bounds for the norm of the basis elements.
Therefore, we appeal to the standard bound.

\begin{definition} Let $\mathcal{K}$ be a subset of $\mathbb{R}^m$ and $\Lambda \subset \mathbb{R}^m$ be a lattice. Then
the {\em successive minima\/} of $\mathcal{K}$ w.r.t. $\Lambda$ are defined for $1\leq k \leq m$ by
\begin{equation*}
  \lambda_k(\mathcal{K},\Lambda):=\inf \left\{\lambda >0 : \lambda \mathcal{K}\ \textnormal{contains}\ k\
   \textnormal{linearly independent points of}\ \Lambda\right\}.
\end{equation*}
\end{definition}

In the following theorem, the term {\em body} denotes a set $\mathcal{K}\subset \mathbb{R}^m$ with non-empty interior such that
$\mathcal{K}$ is contained in the closure of its interior.

\begin{theorem}[Minkowski's Second Theorem]\label{thm:mink2}
  If $\Lambda$ is an $m$-dimensional lattice in $\mathbb{R}^m$ and $\mathcal{K}\subset \mathbb{R}^m$
  is a bounded zero-symmetric convex body with volume $V(\mathcal{K})$, then
  \[
    \lambda_1(\mathcal{K},\Lambda) \dotsm \lambda_m(\mathcal{K},\Lambda) V(\mathcal{K}) \leq 2^m \bm{d}(\Lambda).
  \]
\end{theorem}
\begin{proof}
  See Gruber and Lekkerkerker's monograph~\cite[Theorem~2.16.3]{Gr87}.
\end{proof}

From this theorem we will deduce that $L_g(a_1,a_2)$ must contain either two `short' linearly
independent vectors or one `very short' nonzero vector.
If the first case occurs, we will apply the following result of Bender~\cite{Be62}.

\begin{lemma}\label{le:bender}
  Let $\{\bm{w}_1,\bm{w}_2\}$ be a basis of a lattice $\Lambda\subset\mathbb{R}^2$, and let $0<\vartheta<\pi$
  be the angle between $\bm{w}_1$ and $\bm{w}_2$. Suppose further that
  $\mathcal{C}\subset\mathbb{R}^2$ is a bounded convex set such that the quotient of its area and its perimeter
  is greater than
  \[
    \tfrac{1}{2}\max\left(\|\bm{w}_1\|_2,\|\bm{w}_2\|_2 \sin\vartheta\right).
  \]
  Then $\mathcal{C}$ contains a point of $\Lambda$.
\end{lemma}

For the second case, where we find one vector of `very small' norm in $L_g(a_1,a_2)$, we could
not find an applicable result in the literature that would ensure a lattice point
in the square, so we provide one now.

\begin{figure}[h]
  \centering
  \includegraphics[scale=.6]{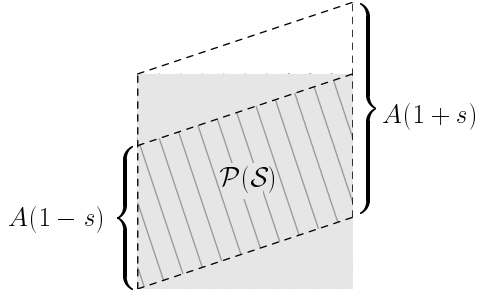}
  \caption{The square $\mathcal{S}$ (shaded) and the parallelogram $\mathcal{P}(\mathcal{S})$ (hatched), which
  lies between two lines of slope $s$ that go through the upper right and the lower left corner of $\mathcal{S}$,
  respectively.}\label{fi:slope}
\end{figure}

\begin{lemma}\label{le:one short vector}
  Let $\Lambda\subset\mathbb{R}^2$ be a lattice and $\bm{r}=(r_1,r_2)$ be a point of $\Lambda$ with
  $\gcd(r_1,r_2)=1$ and $0<r_2\leq r_1$.
  Let further $\mathcal{Q}$ be an open square with sides parallel to the axes and side length $A>0$.
  If $\mathcal{Q}$ contains no point of $\Lambda$, then
  \begin{equation*}
    A \leq \max\left( r_1, \tfrac{\bm{d}(\Lambda)+2r_1r_2}{r_1+r_2}\right).
  \end{equation*}
\end{lemma}
\begin{proof}
  There is a family $\mathfrak{L}$ of parallel equidistant
  lines with slope $s:=r_2/r_1$ such that $\Lambda\subset \bigcup\mathfrak{L}$ and
  the perpendicular distance between two adjacent lines of $\mathfrak{L}$ is $\bm{d}(\Lambda)/\|\bm{r}\|_2$
  \cite[Lemma~III.5]{Ca59}.
  Then the vertical distance between two adjacent lines is $D:=\bm{d}(\Lambda)/r_1$.
  We claim
   \begin{multline}\label{eq:min max}
    \min_{\bm{c}\in\mathbb{R}^2} \max_{\ell\in\mathfrak{L}}\ \left(\textnormal{horizontal length of}\ 
    \ell\cap\mathcal{S}_{A/2}(\bm{c})\right) \\
    =\begin{cases}
      A, & D\leq A(1-s); \\
      \frac{A(1+s)-D}{2s}, & A(1-s) \leq D \leq A(1+s); \\
      0, & D \geq A(1+s).
    \end{cases}
  \end{multline}
  If $D\leq A(1-s)$, then for each square $\mathcal{S}=\mathcal{S}_{A/2}(\bm{c})$ there is a line in $\mathfrak{L}$ that goes
  through the left and the right edge of the square (see Figure~\ref{fi:slope}). This settles the first
  case in the right hand side of \eqref{eq:min max}.

  If $D$ is larger than $A(1+s)$, there
  is a square that is not intersected by any line from $\mathfrak{L}$.

  We are left with the intermediate case
  $A(1-s) \leq D \leq A(1+s)$. To achieve the minimum in~\eqref{eq:min max}, we must certainly place $\mathcal{S}$
  such that there is no line from $\mathfrak{L}$ in the parallelogram $\mathcal{P}(\mathcal{S})$ of Figure~\ref{fi:slope}.
  But then there is always a line $\ell\in\mathfrak{L}$ that intersects $\mathcal{S}\backslash \mathcal{P}(\mathcal{S})$, say in the
  upper triangle of $\mathcal{S}\backslash \mathcal{P}(\mathcal{S})$. If no line intersects the lower triangle of
  $\mathcal{S}\backslash \mathcal{P}(\mathcal{S})$, we can make the maximum in~\eqref{eq:min max} smaller by pushing
  $\mathcal{S}$ downwards. The smallest possible value of the maximum is achieved as soon as the intersections of $\mathcal{S}$ with
  $\ell$ and the line from $\mathfrak{L}$ just below $\ell$ have equal length. It is easy to see that these intersections both
  have horizontal length $(A(1+s)-D)/2s$.

  Now that \eqref{eq:min max} is established, let
  $\mathcal{Q}$ be an open square with sides parallel to the axes and side length
  \begin{equation}\label{eq:estimateA}
    A > \max\left( r_1, \tfrac{\bm{d}(\Lambda)+2r_1r_2}{r_1+r_2}\right).
  \end{equation}
  Our goal is to show $\mathcal{Q}\cap \Lambda\neq \emptyset$.
  If the first case in the right hand side of~\eqref{eq:min max} occurs, we are well off:
  Since $A > r_1$, the line segment in $\mathcal{Q}\cap\bigcup\mathfrak{L}$
  of horizontal length $A$ must contain a point of $\Lambda$.  
  The third case in~\eqref{eq:min max} cannot happen, since it would imply $\bm{d}(\Lambda)\geq A(r_1+r_2)$,
  contradicting \eqref{eq:estimateA}.
  As for the second case, $A > \tfrac{\bm{d}(\Lambda)+2r_1r_2}{r_1+r_2}$ implies
  \[
    r_1 < \frac{A(r_1+r_2)-\bm{d}(\Lambda)}{2r_2} = \frac{A(1+s)-D}{2s},
  \]
  hence $\mathcal{Q}\cap \Lambda\neq \emptyset$.
\end{proof}

\begin{proof}[Proof of Proposition~\ref{pr:g}]
  We begin this proof, which is the core of the proof of Theorem~\ref{thm:main},
  by settling the cases where $g$ is at most $9$. The only numbers to consider are $g=7,8,9$,
  since for smaller $g\neq 5$ there are no $a_1,a_2$ that satisfy the requirements of Theorem~\ref{thm:genrat}
  (and hence of Proposition~\ref{pr:g}).
  First let $g=7$. If we have proved the desired result for a pair $(a_1,a_2)$, we need not
  consider the five pairs
  \[
    (a_2,a_1), (g-a_1,a_2), (a_1,g-a_2),
    (g-a_2,a_1) \quad \text{and} \quad (a_2,g-a_1)
  \]
  any more by Lemma~\ref{le:symm}. It is readily seen that
  under our restrictions on $a_1$, $a_2$ all lattices $L_7(a_1,a_2)$
  are equal to
  $L_7(1,3)$ modulo these symmetries. Similarly, for $g=8$ and $g=9$
  it suffices to consider $L_8(3,1)$ and $L_9(2,1)$, respectively.
  In all three cases it is easy to verify the desired result.

  From now on we assume $g\geq 10$.
  Put $\Lambda:=L_g(a_1,a_2)$ and let
  \[
    \mathcal{B} := \left\{ \bm{x} \in\mathbb{R}^2 : \|\bm{x}\|_2 \leq 1\right\}
  \]
  be the unit circle.
  It is not difficult to see \cite[Section~2]{Le94} that the determinant of $\Lambda$ is $\bm{d}(\Lambda)=g$.
  Then Theorem~\ref{thm:mink2} shows
  \[
    \lambda_1(\mathcal{B},\Lambda)\lambda_2(\mathcal{B},\Lambda)\pi \leq 4g.
  \]
  First suppose $\lambda_2(\mathcal{B},\Lambda)<g/4$.
  The quotient of the area of $\mathcal{S}_{g/4}(\bm{c})$
  and its perimeter is $\tfrac{g^2}{4}/2g=g/8$, hence we can apply Lemma~\ref{le:bender}.

  If, on the other hand, $\lambda_2(\mathcal{B},\Lambda)\geq g/4$,
  then we have $\lambda_1(\mathcal{B},\Lambda)\leq \tfrac{16}{\pi}$, which provides us
  with a nonzero point $\bm{r}\in\Lambda$ with $\|\bm{r}\|_2\leq \tfrac{16}{\pi}$.
  W.l.o.g.\ assume that $\bm{r}$ satisfies $\gcd(r_1,r_2)=1$ and $0<r_2\leq r_1$.
  According to Lemma~\ref{le:one short vector}, it suffices to show
  \[
    \frac{g}{2} > \frac{g+2r_1r_2}{r_1+r_2},
  \]
  i.e.
  \begin{equation}\label{eq:r1 r2}
    4r_1r_2 < g(r_1+r_2-2).
  \end{equation}
  This inequality is satisfied for $g\geq 10$ and
  \[
    \bm{r} \in \left\{ (2,1),(3,1),(4,1),(2,2),(3,2),(4,2),(3,3),(4,3)\right\},
  \]
  which are all values of $\bm{r}$ in question. Observe that $a_1\not\equiv a_2\pmod{g}$ implies
  $(1,1)\notin\Lambda$.
\end{proof}

This completes the proof of Theorems~\ref{thm:genrat} and~\ref{thm:xi0,xi1}.
We remark that the successive minima approach from the preceding proof can be applied to the case
of distinct denominators $b_1$, $b_2$, too. However, the number of special
cases that have to be checked separately is much larger than for equal denominators.

\section{Completion of the proof of the main theorem}\label{se:compl}

If $\xi_2$ from Theorem~\ref{thm:xi0,xi1} equals $\tfrac{1}{2}$, which corresponds to $\arg{\alpha_2}=\theta_2=\pi$
and thus to a real negative dominating characteristic root $\alpha_2$ in
Theorem~\ref{thm:main}, then
the squares centred at $\bm{c}=(c_1,\tfrac{1}{4})$ or $\bm{c}=(c_1,\tfrac{3}{4})$, $c_1\in\mathbb{R}$, do not contain
any point $n\bm{\xi}\bmod 1$.
But in this case we need not consider all squares:
\[
  w_2 \sin(n \pi+\varphi_2) = (-1)^n w_2 \sin{\varphi_2},
\]
hence $\varphi_2$ can be absorbed in $w_2$, and
we retain full generality if we assign a convenient value to $\varphi_2$.

\begin{proposition}\label{pr:eta real}
  Let $\xi_1 \in \left]0,1\right[\backslash \{\tfrac{1}{2}\}$. Then
  for all $c_1 \in \mathbb{R}$ there are $\epsilon > 0$ and  $c_2 \in \mathbb{R}$ such that
  for infinitely many $n\in\mathbb{N}$
  \[
    n(\xi_1,\tfrac{1}{2})\bmod 1\ \in\  \mathcal{S}_{1/4-\epsilon}(\bm{c}) \bmod 1.
  \]
\end{proposition}
\begin{proof}
  Let $c_1$ be an arbitrary real number with $0\leq c_1<1$.
  It suffices to show that there is $n\in\mathbb{N}$ with
  \[
    n\xi_1\bmod 1\ \in\ \left]c_1 - \tfrac{1}{4}, c_1 + \tfrac{1}{4}\right[ \bmod 1.
  \]
  If $\xi_1$ is irrational, this follows immediately from Theorem~\ref{thm:kronecker}.
  Now let $\xi_1$ be a rational number with denominator $b_1 > 0$.
  The desired result follows from $1/b_1 < \tfrac{1}{2}$ and
  \[
    \left\{ n\xi_1\bmod 1 : n\in\mathbb{N} \right\} = \left\{ 0,\tfrac{1}{b_1},\dots,\tfrac{b_1-1}{b_1}\right\}.
  \]
\end{proof}

Applying Proposition~\ref{pr:eta real} with $\xi_1=\theta_1/2\pi$, where $\theta_k=\arg{\alpha_k}$ as usual,
settles the case of Theorem~\ref{thm:main2} where $\alpha_2$ is a negative
real number. Clearly, the same argument applies if $\alpha_1$ is negative real and $\alpha_2$ is complex.

Finally let us see what happens if the pair of denominators of $(\theta_1/2 \pi,\theta_2/2 \pi)$ in
Theorem~\ref{thm:main2} is $(5,5)$, $(6,3)$ or $(8,4)$.
According to Theorem~\ref{thm:genrat}, our argument with lattice points in squares fails for these
values. Fortunately, it is straightforward to show directly that the purely periodic
sequences $(b(n))$ arising from these values oscillate.

Once again we can appeal to the symmetries noted in Lemma~\ref{le:symm}.
Indeed, swapping $\theta_1$ and $\theta_2$ does no harm,
and the sign of $\theta_k$ can be absorbed in $w_1$ and $\varphi_1$.
It turns out that for each of the three pairs of denominators it suffices to consider one pair of numerators,
namely $(a_1,a_2)=(2,2)$ for the denominators $(6,3)$ and $(8,4)$,
and $(a_1,a_2)=(4,2)$ for the denominators $(5,5)$.
In the following proposition, the cases (i), (ii) and (iii) correspond
to the pairs of denominators $(6,3)$, $(8,4)$ and $(5,5)$, respectively.
 
\begin{proposition}
  Let $w_1$, $w_2$ be nonzero real numbers and $\varphi_1$, $\varphi_2$ be real numbers. Then the sequence
  defined by
  \[
    b(n) := w_1 \sin(n\theta_1+\varphi_1) + w_2 \sin(n\theta_2+\varphi_2), \quad n\geq 0,
  \]
  has a positive and a negative entry for each of the following values of $(\theta_1,\theta_2)$.
  \begin{flalign*}
    (i) & \quad \theta_1=\tfrac{1}{3}\pi, \theta_2=\tfrac{2}{3}\pi & \\
    (ii) & \quad \theta_1=\tfrac{1}{4}\pi, \theta_2=\tfrac{1}{2}\pi & \\
    (iii) & \quad \theta_1=\tfrac{4}{5}\pi, \theta_2=\tfrac{2}{5}\pi &
  \end{flalign*}
\end{proposition}
\begin{proof}
  It suffices to consider $w_1=1$.
  We set $s_i:=\sin\varphi_i$ and $c_i:=\cos\varphi_i$ for $i=1,2$,
  suppose $b(n)\geq 0$ for all $n\geq 0$ and derive a contradiction.

  (i) From
  \begin{align*}
    b(0)+b(1)+b(5) &=2s_1, \\
    b(2)+b(3)+b(4) &=-2s_1, \\
    b(0)+b(3) &=2s_2w_2, \\
    b(1)+b(2)+b(4)+b(5) &=-2s_2w_2
  \end{align*}
  we deduce $s_1=s_2w_2=0$, so $s_2=0$ or $w_2=0$. If $w_2=0$, then
  \begin{align*}
    b(1)& = \tfrac{1}{2}\sqrt{3}c_1, \\
    b(4)& = -\tfrac{1}{2}\sqrt{3}c_1
  \end{align*}
  implies $c_1=0$, which contradicts $c_1^2+s_1^2=1$.
  If, on the other hand,  $s_2=0$, then $c_1=0$ follows from
  \begin{align*}
    b(1)+b(2)=\sqrt{3}c_1, \\
    b(4)+b(5)=-\sqrt{3}c_1.
  \end{align*}

  (ii) Now
  \begin{align*}
    b(1)+b(3) &=\sqrt{2}c_1, \\
    b(5)+b(7) &=-\sqrt{2}c_1, \\
    b(1)+b(7) &=\sqrt{2}s_1, \\
    b(3)+b(5) &=-\sqrt{2}s_1
  \end{align*}
  shows $c_1=s_1=0$, again contradicting $c_1^2+s_1^2=1$.

  (iii) Since
  \begin{align*}
    b(0) &=s_1+s_2w_2, \\
    b(1)+b(2)+b(3)+b(4) &=-s_1-s_2w_2,
  \end{align*}
  we have $s_2w_2=-s_1$. Then we obtain
  \begin{align*}
    b(1)+b(4) &=-\sqrt{5}s_1, \\
    b(2)+b(3) &=\sqrt{5}s_1,
  \end{align*}
  hence $s_1=0$ and $c_1=\pm 1$. Therefore $w_2=0$ or $s_2=0$.
  If $w_2=0$, then the values
  \[
    b(1) = \tfrac{1}{2} \sqrt{\tfrac{1}{2}(5-\sqrt{5})} c_1
  \]
  and
  \[
    b(2) = -\tfrac{1}{2} \sqrt{\tfrac{1}{2}(5+\sqrt{5})} c_1
  \]
  have opposite signs. It remains to consider the case $s_2=0$.
  For each of the four possible values of $(c_1,c_2)=(\pm 1,\pm 1)$ the inequalities
  $b(0)\geq 0,\dots,b(4)\geq 0$ form a linear system of inequalities in $w_2$. It is easy
  to check that none of these four systems are solvable.
\end{proof}

The proof of Theorem~\ref{thm:main2} is complete, hence Theorem~\ref{thm:main} is established.

\section{A positive real root}\label{se:pos real}

If one of the dominating characteristic roots $\alpha_1,\dots,\alpha_m$ is real positive, Conjecture~\ref{cj:equ abs}
is not applicable. Consider the sequence defined by
\begin{equation}\label{eq:a pos real}
  a(n) := \sum_{k=1}^m w_k\sin(n\theta_k+\varphi_k)+1+o(1), \quad n\geq 0,
\end{equation}
where $\theta_1,\dots,\theta_m$, $w_1,\dots,w_m$ are nonzero real numbers and $\varphi_1,\dots,\varphi_m$ are real numbers.
Here and throughout this section we assume that
the coefficient of the real positive root is positive (and
thus w.l.o.g.\ equals one). Analogous
considerations apply for a negative coefficient.
The behaviour of $(a(n))$
depends on how $1$ compares to
\begin{equation}\label{eq:s real root}
  S:=-\inf_{n\geq 0} \ \sum_{k=1}^m w_k\sin(n\theta_k+\varphi_k)\ \in\ \ ]-W,W],
\end{equation}
where $W=\sum_{k=1}^m|w_k|$.
The sequence $(a(n))$ is positive for large $n$ if $S<1$ (in particular, if $W<1$),
and it oscillates if $S>1$.
If $S=1$, the behaviour of $(a(n))$ depends on how well
$\sum_{k=1}^m w_k\sin(n\theta_k+\varphi_k)$ approximates $-1$ and possibly on the $o(1)$ term.

The preceding discussion gives a handy criterion only for $W<1$, which was already
obtained by Burke and Webb \cite{Bu81}.
For $W\geq 1$ we confine ourselves to showing
how $(a(n))$ behaves for almost all values of the parameters $\theta_k$, $\varphi_k$ and $w_k$.

\begin{lemma}\label{le:almost xi}
  Let $\bm{\alpha}\in\mathbb{R}^m$ and let $(\psi(n))$ be a sequence
  of positive real numbers such that $\sum_{n\geq 0} \psi(n)^m$ converges.
  Then the set of inequalities
  \[
    (n\xi_k-\alpha_k)\bmod 1 < \psi(n), \quad 1\leq k\leq m
  \]
  has infinitely many solutions $n\in\mathbb{N}$ for almost no $\bm{\xi}\in\mathbb{R}^m$.
\end{lemma}
\begin{proof}
  See Cassels \cite[Lemma~VII.2.1]{Ca57}.
\end{proof}

In order to apply the following theorem we require the dominating characteristic roots to be simple.
This assumption makes the remainder term $r(n)$ go to zero exponentially. Parts (i) and (iii)
hold for multiple roots as well, since they only require $r(n)=o(1)$.
Our proof of part (ii), however, breaks down for $m=1$ in case of a multiple root,
because then we can ensure only  $r(n)=O(n^{-1})$ and this leads to a divergent series in Lemma~\ref{le:almost xi}.

\begin{theorem}\label{thm:almost no theta}
  Let $w_1,\dots,w_m$ be nonzero real numbers with $W:=\sum_{k=1}^m|w_k|$,
  $\varphi_1,\dots,\varphi_m$ be real numbers and $(r(n))$ be a real sequence
  with $r(n)=O(\omega^n)$ for some $0<\omega<1$.
  \begin{itemize}
    \item[(i)] If $W<1$, then for all $\bm{\theta}\in\mathbb{R}^m$ the sequence $(a(n))$ defined by
    \[
      a(n) := \sum_{k=1}^m w_k\sin(n\theta_k+\varphi_k)+1+r(n)
    \]
    is positive for large $n$.
    \item[(ii)] If $W=1$, then for almost all $\bm{\theta}\in\mathbb{R}^m$ the sequence $(a(n))$ is positive for large $n$.  
    \item[(iii)] If $W>1$, then $(a(n))$ oscillates for almost all $\bm{\theta}\in\mathbb{R}^m$.
  \end{itemize}
\end{theorem}
\begin{proof}
  (i) is clear.
  (iii) follows from the $m$-dimensional version of Theorem~\ref{thm:kronecker} (see Section~\ref{se:conclusion}), because $\theta_1/2\pi,\dots,\theta_m/2\pi,1$ are linearly
  independent over the rationals for almost all $\bm{\theta}$. We proceed to prove (ii).
  Suppose $a(n)\leq 0$ for all $n$ in an infinite set $I\subseteq \mathbb{N}$. To make $a(n)$ non-positive, $\sin(n\theta_k+\varphi_k)$
  has to be very close to $-1$ if $w_k>0$ and very close to $1$ if $w_k<0$. To be precise, we must have
  \[
    \lim_{\substack{n\to\infty \\ n\in I}}f(n) = 0
  \]
  for
  \[
    f(n) := (f_1(n),\dots,f_m(n))
  \]
  with
  \[
    f_k(n) :=
    \begin{cases}
      (n\theta_k +\varphi_k - \tfrac{1}{2}\pi)\bmod 2\pi, & w_k < 0; \\
      (n\theta_k + \varphi_k - \tfrac{3}{2}\pi)\bmod 2\pi, & w_k > 0.
    \end{cases}
  \]
  By Taylor expansion, we obtain
  \begin{align*}
     \sum_{k=1}^m w_k\sin(n\theta_k+\varphi_k)+1
    &= -\sum_{k=1}^m|w_k| + \frac{1}{2}\sum_{k=1}^m|w_k|f_k(n)^2 + O(\sum_{k=1}^m f_k(n)^4) + 1 \\
    &= \frac{1}{2}\sum_{k=1}^m|w_k|f_k(n)^2 + O(\sum_{k=1}^m f_k(n)^4) \quad \text{as}\ n\to\infty\ \text{in}\ I. 
  \end{align*}
  Removing finitely many elements from $I$ if necessary, we thus have
  \[
    \sum_{k=1}^m w_k\sin(n\theta_k+\varphi_k)+1 > \frac{w}{3}\sum_{k=1}^m f_k(n)^2, \quad n\in I,
  \]
  where $w:=\min_{1\leq k\leq m}|w_k|>0$. Since $a(n)\leq 0$ for $n\in I$, this implies
  \[
    \sum_{k=1}^m f_k(n)^2 < -\tfrac{3}{w}r(n) = O(\omega^n), \quad n\in I,
  \]
  hence for $1\leq k\leq m$
  \[
    f_k(n) = O(\omega^{n/2})\quad  \text{as} \ n\to\infty \ \text{in}\ I.
  \]
  According to Lemma~\ref{le:almost xi} this holds for almost no $\bm{\theta}$.
\end{proof}

Finer questions may be asked about the sets of measure zero alluded to in Theorem~\ref{thm:almost no theta}.
As for part (ii) of the theorem, we note that there are $\varphi_1,\dots,\varphi_m$,
$r(n)$ and infinitely many $\bm{\theta}$ such that $(a(n))$ oscillates for all nonzero $w_1,\dots,w_m$ with $W=1$.
To see this, define
\[
  \varphi_k :=
  \begin{cases}
    \tfrac{1}{2}\pi, & w_k < 0; \\
    \tfrac{3}{2}\pi, & w_k > 0,
  \end{cases}
\]
let $\bm{\theta}\in\mathbb{Q}^m$ be arbitrary
and $r(n):=(-\omega)^{n+1}$ for some $0<\omega<1$. Then $a(n)\geq (-\omega)^{n+1}=\omega^{n+1}>0$
for odd $n$, and $a(n)=-\omega<0$ if $n$ is two times a common multiple of the denominators of $\theta_1,\dots,\theta_m$.

The preceding example is a special case of the following proposition, which completely
describes the behaviour of $(a(n))$ under the assumptions of part (ii) of Theorem~\ref{thm:almost no theta}
and the additional constraint $\theta_k/\pi\in\mathbb{Q}$ for $1\leq k\leq m$.

\begin{proposition}\label{pr:theta rat}
  Let $\theta_1,\dots,\theta_m$ be real numbers such that $\theta_k/2\pi$
  is a rational number $a_k/b_k$ for $1\leq k\leq m$, let $\varphi_1,\dots,\varphi_m$ be real numbers,
  let $w_1,\dots,w_m$ be nonzero real numbers with $\sum_{k=1}^m |w_k|=1$ and define
  \[
    a(n) := \sum_{k=1}^m w_k\sin(n\theta_k+\varphi_k)+1 + o(1), \quad n\geq 0.
  \]
  \begin{itemize} 
    \item[(i)] If there is a $k$ such that $\varphi_k/\pi$ is irrational, then $(a(n))$ is positive for large $n$.
    \item[(ii)] Suppose that $\varphi_k/2\pi$ is a rational number $c_k/d_k$ for $1\leq k\leq m$.
    If for all $1\leq k,l\leq m$
    \begin{equation}\label{eq:theta rat congr}
      b_k(A_kd_k - 4c_k) \equiv b_l(A_ld_l - 4c_l) \pmod{4\gcd(d_kb_k,d_lb_l)}
    \end{equation}
    with
    \[
      A_k :=
      \begin{cases}
        1, & w_k < 0; \\
        3, & w_k > 0,
      \end{cases}
    \]
    then there are infinitely many $n$ with $b(n)=0$, where
    \[
      b(n) := \sum_{k=1}^m w_k\sin(n\theta_k+\varphi_k) + 1 \geq 0,
    \]
    and the behaviour of $(a(n))$ depends in an obvious way on the sign of the $o(1)$ term
    for these $n$.
    If there are $k$, $l$ such that \eqref{eq:theta rat congr} does not hold, then $(a(n))$ is positive for large $n$.
  \end{itemize}
\end{proposition}
\begin{proof}
  The purely periodic sequence $(b(n))$ satisfies $b(n)\geq 0$ for all $n\geq 0$. If none of its finitely many values are zero, then
  $(a(n))$ is positive for large $n$. We have $b(n)=0$ if and only if
  $\sin(n\theta_k+\varphi_k)$ equals $1$ for the $k$'s with $w_k<0$ and
  $-1$ for the $k$'s with $w_k>0$, i.e.
  \[
    n\theta_k + \varphi_k \equiv \tfrac{1}{2}A_k\pi \pmod {2\pi}, \quad 1\leq k\leq m,
  \]
  which is equivalent to
  \[
    n\tfrac{a_k}{b_k} + \tfrac{\varphi_k}{2\pi} \equiv \tfrac{1}{4}A_k \pmod{1}, \quad 1\leq k\leq m.
  \]
  Clearly, this cannot hold if one of the $\varphi_k/\pi$ is irrational.
  Under the assumption of part (ii), we are lead to the system of congruences
  \[
    4d_ka_k n \equiv b_k(A_kd_k - 4c_k) \pmod{4d_kb_k}, \quad 1\leq k\leq m.
  \]
  Now the result follows from Theorem~\ref{thm:gcrt}.
\end{proof}

\section{Conclusion and future directions}\label{se:conclusion}

In order to extend our approach to Conjecture~\ref{cj:equ abs} to $m$ dominating characteristic roots,
we need to show that infinitely many $n(\xi_1,\dots,\xi_m)\bmod 1$ lie in any given $m$-dimensional hypercube
(modulo one) with side length $\tfrac{1}{2}-\epsilon$.
Theorem~\ref{thm:kronecker} generalizes in the following way \cite[Theorem~III.5.IV]{Ca57}:
The points $n\bm{\xi}\bmod 1$ lie dense in the set of all
$\bm{x}\in\mathbb{R}^m$ that satisfy $\langle \bm{u},\bm{x} \rangle\in\mathbb{Z}$ for all
integer vectors $\bm{u}$ with $\langle \bm{u},\bm{\xi} \rangle\in\mathbb{Z}$.
Again the case of rational $\xi_1=a_1/b_1,\dots,\xi_m=a_m/b_m$ with equal denominators
$b_1=\dots=b_m=g$ will be the crux of the proof.
This case seems to become more and more difficult for fixed denominator $g$ as $m$ increases, since the set
\begin{equation}\label{eq:set m}
  \left\{n(\tfrac{a_1}{g},\dots,\tfrac{a_m}{g})\bmod 1 : n\in\mathbb{N}\right\}
\end{equation}
has $g$ elements for all $m$, whereas the volume of the hypercube
is $(\tfrac{1}{2}-\epsilon)^m$.
Theorem~\ref{thm:mink2} is certainly a valuable tool.
Hadwiger \cite{Ha70} has extended Bender's two-dimensional result (Lemma~\ref{le:bender})
that we used in the proof of Proposition~\ref{pr:g}
to arbitrary dimension $m$. A significant extension of Lemma~\ref{le:one short vector} is still needed.
Anyway it is conceivable that the exceptional rational values $a_1/b_1,\dots,a_m/b_m$
for which the hypercube might contain no point of \eqref{eq:set m} become unmanageable
as $m$ increases (Cf. Section~\ref{se:compl}).

Our results on a positive real characteristic root leave ample room for refinement. For instance,
one could try to extend part (ii) of Theorem~\ref{thm:almost no theta} to the case of a multiple root
or to continue the discussion begun in Proposition~\ref{pr:theta rat} by relaxing the
requirement that all $\theta_k/\pi$ be rational.

The Skolem--Mahler--Lech Theorem \cite{Ev03} describes
the structure of the zero set $\{n:a(n)=0\}$ of a recurrence sequence.
It is the union of a finite set and finitely many arithmetic progressions.
There might be an analogue of this result for the set $\{n:a(n)>0\}$.

Finally, we have excluded algorithmics so far. We do not know whether the positivity of recurrence sequences
is a decidable problem. Proving Conjecture~\ref{cj:equ abs} and
giving an effective criterion instead of the metric Theorem~\ref{thm:almost no theta}
would lead to a decision procedure.


\begin{thebibliography}{00}

\bibitem{An00}
Andrews, G. E., Askey, R. and Roy, R., 2000,
{\itshape Special Functions.}
Cambridge University Press.

\bibitem{Ze90}
Zeilberger, D., 1990,
A holonomic systems approach to special functions identities.
{\itshape J.\ Comput.\ Appl.\ Math.\/} {\bfseries 32}, 321--368.

\bibitem{Ev03}
Everest, G., van der Poorten, A., Shparlinski, I. and Ward, T., 2003,
{\itshape Recurrence Sequences},
American Mathematical Society, Providence RI.

\bibitem{Bu81}
Burke, J. R. and Webb, W. A., 1981,
Asymptotic behavior of linear recurrences.
{\itshape Fibonacci Quart.\/} {\bfseries 19}, 318--321.

\bibitem{Na90}
Nagasaka, K. and Shiue, J., 1990,
Asymptotic positiveness of linear recurrence sequences.
{\itshape Fibonacci Quart.\/} {\bfseries 28}, no.\ 4, 340--346.

\bibitem{Ca59}
Cassels, J. W. S., 1959,
{\itshape An Introduction to the Geometry of Numbers},
Springer Berlin, Goettingen, Heidelberg.

\bibitem{Kn94}
Knuth, D. E., 1994,
{\itshape Concrete Mathematics\/} (2nd edn),
Addison Wesley, Reading MA.

\bibitem{Ni63}
Niven, I., 1963,
{\itshape Diophantine Approximations\/},
Interscience, New York.

\bibitem{Ro97}
Rote, G., 1997,
Finding a shortest vector in a two-dimensional lattice modulo $m$.
{\itshape Theoret.\ Comput.\ Sci.\/} {\bfseries 172}, 303--308.

\bibitem{Kn98}
Knuth, D. E., 1998,
{\itshape The Art of Computer Programming, vol.\ 2, Seminumerical Algorithms\/} (3rd edn),
Addison Wesley, Reading MA.

\bibitem{Le94}
Lempel, M., and Paz, A., 1994,
An algorithm for finding a shortest vector in a two-dimensional modular lattice.
{\itshape Theoret.\ Comput.\ Sci.\/} {\bfseries 125}, 229--241.

\bibitem{Gr87}
Gruber, P. M. and Lekkerkerker, C. G., 1987,
{\itshape Geometry of Numbers} (2nd edn),
North-Holland, Amsterdam.

\bibitem{Be62}
Bender, E. A., 1962,
Area-perimenter relations for two-dimensional lattices.
{\itshape Amer.\ Math.\ Monthly} {\bfseries 69}, 742--744.

\bibitem{Ca57}
Cassels, J. W. S., 1957,
{\itshape An Introduction to Diophantine Approximation},
Cambridge University Press.

\bibitem{Ha70}
Hadwiger, H., 1970,
Volumen und Oberfl\"{a}che eines Eik\"{o}rpers, der keine Gitterpunkte \"{u}berdeckt.
{\itshape Math.\ Z.\/} {\bfseries 116}, 191--196. 


\end{thebibliography}
\end{document}